\documentclass[10pt]{article}
\usepackage{amsfonts}
\usepackage{mathrsfs}
\usepackage{mathrsfs}
\usepackage{amsmath}
\usepackage{amsfonts}
\usepackage{amsfonts}
\usepackage{amsfonts}
\usepackage{amsmath}
\usepackage{amssymb}
\usepackage{amsfonts}
\usepackage{amsthm}
\usepackage{mathrsfs}
\usepackage[all]{xy}
\usepackage{graphicx}
\newtheorem{theorem}{Theorem}[section]

\newtheorem{lemma}[theorem]{Lemma}

\theoremstyle{definition}

\begin{document}

\title{A Bourgain type bilinear estimate for a class of water-wave models }
% Enter your title between curly braces
\author{Qifan Li}
\date{September,\ \ 2010}          % Enter your date or \today between curly braces
\maketitle We consider the general form of the equation of
water-wave models on torus $\mathbb{T}$
\begin{equation}\label{1}\partial_tu+\sum_{k=1}^Nb_k\partial_x^{2k+1}u+
Q(u,\partial_xu,\cdots\partial_x^{2N+1}u)=0\end{equation} where $Q$
denotes nonlinear term of the equation and $b_k$s are real
constants. This equation was first introduced and studied by
\cite{l1} on the real line. The symbol of the linear partial
differential operator $\mathcal{L}=\sum_{k=1}^Nb_k\partial_x^{2k+1}$
is
$$m(\xi)=\sum_{k=1}^{N}b_k(2\pi
i\xi)^{2k+1}=P(\xi)i$$ where $P(\xi)=\sum_{k=1}^{N}c_k\xi^{2k+1}$
and $c_k=(-1)^k(2\pi)^{2k+1}b_k.$ We assume $c_k\leq0$ for all the
$k\geq1$. In particular, generalized Kawahara equations \cite{B1}
and fifth-order KdV equations \cite{K} are the special cases
satisfying this condition. We introduce the Bourgain space-time
space $X^{s,b}$ with the norm
\begin{equation}\label{2}\|u\|_{X^{s,b}}^2=\sum_{n\in\mathbb{Z}}\int_{\mathbb{R}}(1+|\lambda+P(n)|)^{2b}(1+|n|)^{2s}|
\hat{u}(n,\lambda)|^2d\lambda.\end{equation}

We are going to establish a bilinear estimate which extends the
proposition 7.15 in \cite{Bo} to the higher order derivatives.

\begin{lemma}Let functions $u,\ v:\ \mathbb{T}
         \times[0,T]\rightarrow\mathbb{R}$, we
          have the bilinear estimate\begin{equation}\label{3}\|uv\|_{L^2_xL^2_t}\lesssim\|u\|_{X^{0,\frac{N+1}{4N+2}}}
\|v\|_{X^{0,\frac{N+1}{4N+2}}}.\end{equation}\end{lemma}

\begin{proof}Let \[\begin{split}\label{5}A_m&=\left(\int_{\mathbb{R}}(1+|\xi|)^{\frac{N+1}{2N+1}}
|\widehat{u_m}(\xi)|^2d\xi
\right)^{1/2}\\B_m&=\left(\int_{\mathbb{R}}(1+|\xi|)^{\frac{N+1}{2N+1}}|\widehat{v_m}(\xi)|^2d\xi
\right)^{1/2}\end{split}\]and\begin{equation}\label{5}u(x,t)=\sum_{m\in\mathbb{Z}}e^{2\pi
i(mx-P(m)t)}u_m(t),\ \ \ v(x,t)=\sum_{m\in\mathbb{Z}}e^{2\pi
i(mx-P(m)t)}v_m(t).\end{equation} The proof of (\ref{3}) is reduced
to show $\|uv\|_{L^2_xL^2_t}\lesssim\|A_m\|_{l_n^2}\|B_m\|_{l_n^2}$.

Let quadratic polynomial $Q(m,l)=P(m+l)-P(m)-P(l)$. We have
$$(u\bar{v})(x,t)=\sum_{l\in\mathbb{Z}}e^{2\pi i(lx-P(l)t)}\sum_{m\in\mathbb{Z}}e^{2\pi iQ(m,l)t}
(u_m\overline{v_{m+l}})(t)$$and\begin{equation}\label{5}
\|u\bar{v}\|_{L^2_xL^2_t}^2=\sum_{l\in\mathbb{Z}}\int_{\mathbb{R}}\left|\sum_{m\in\mathbb{Z}}e^{2\pi
iQ(m,l)t} (u_m\overline{v_{m+l}})(t)\right|^2dt.\end{equation}

For a integer $j>0$, define Paley-Littlewood operator
$$\Delta_jf(t)=(\mathbf{1}_{2^{j-1}\leq|\xi|\leq2^{j}}\hat{f}(\xi))^{\vee}(t)$$
where $\mathbf{1}_{2^{j-1}\leq|\xi|\leq2^{j}}$ denotes the
characteristic function on the set
$[2^{j-1},2^j]\cup[-2^j,-2^{j-1}]$, and set
$$\Delta_0f(t)=(\mathbf{1}_{|\xi|\leq1}\hat{f}(\xi))^{\vee}(t).$$
We have the Paley-Littlewood decomposition for $u_m$ and $v_m$
$$u_m=\sum_{p\geq0}\Delta_{p}u_m,\ \ \ \ \ \ \ v_m=\sum_{q\geq0}\Delta_{q}v_m.$$

We assert that (\ref{5}) can be estimated by
\begin{equation}\begin{split}\label{6}
\|u\bar{v}\|_{L^2_xL^2_t}^2&\lesssim\sum_{l\in\mathbb{Z}}\left(\sum_{q\geq
p}\left\|\sum_{m\in\mathbb{Z}}e^{2\pi iQ(m,l)t}
(\Delta_{p}u_m\overline{\Delta_{q}v_{m+l}})(t)\right\|_{L_t^2}\right)^2\\&+\sum_{l\in\mathbb{Z}}\left(\sum_{q\geq
p}\left\|\sum_{m\in\mathbb{Z}}e^{2\pi iQ(m,l)t}
(\Delta_{p}v_m\overline{\Delta_{q}u_{m+l}})(t)\right\|_{L_t^2}\right)^2.\end{split}\end{equation}Since
it is easy to see that
$$\|u\bar{v}\|_{L^2_xL^2_t}^2\leq\sum_{l\in\mathbb{Z}}\left(\sum_{q\in\mathbb{Z}}\sum_{p\in\mathbb{Z}}
\left\|\sum_{m\in\mathbb{Z}}e^{2\pi iQ(m,l)t}
(\Delta_{p}u_m\overline{\Delta_{q}v_{m+l}})(t)\right\|_{L_t^2}\right)^2.$$For
the summation over $q\leq p$ inside the brackets, take
$m^{\prime}=m+l$ and $l^{\prime}=-l$ we have
$Q(m^{\prime}+l^{\prime},-l^{\prime})=P(m^{\prime})-P(m^{\prime}+l^{\prime})-P(-l^{\prime})=
-Q(m^{\prime},l^{\prime})$. By this reason, we can
write\begin{equation}\begin{split}\label{7}\sum_{l\in\mathbb{Z}}\left(\sum_{q\leq
p} \left\|\sum_{m\in\mathbb{Z}}e^{2\pi iQ(m,l)t}
(\Delta_{p}u_m\overline{\Delta_{q}v_{m+l}})(t)\right\|_{L_t^2}\right)^2&=\sum_{l^{\prime}\in\mathbb{Z}}\left(\sum_{q\leq
p} \left\|\sum_{m^{\prime}\in\mathbb{Z}}e^{2\pi
iQ(m^{\prime}+l^{\prime},-l^{\prime})t}
(\Delta_{p}u_{m^{\prime}+l^{\prime}}\overline{\Delta_{q}v_{m^{\prime}}})(t)\right\|_{L_t^2}\right)^2\\&=
\sum_{l^{\prime}\in\mathbb{Z}}\left(\sum_{q\leq p}
\left\|\sum_{m^{\prime}\in\mathbb{Z}}\overline{e^{2\pi
iQ(m^{\prime},l^{\prime})t}
(\overline{\Delta_{p}u_{m^{\prime}+l^{\prime}}}\Delta_{q}v_{m^{\prime}})(t)}\right\|_{L_t^2}\right)^2.
\end{split}\end{equation}This proves the estimate (\ref{6}).

We will need the pointwise
estimate\begin{equation}\label{8}\sum_m|\Delta_iu_m|^2\leq2^j\left(\sum_m\|\Delta_iu_m\|_{L^2_t}^2\right).\end{equation}
To prove (\ref{8}), we shall use Jensen inequality and Plancherel
theorem.
\[\begin{split}\sum_m|\Delta_iu_m|^2&=\sum_m2^{2j}\left|\frac{1}{2^j}\int_{2^{j-1}\leq|\xi|\leq2^j}\widehat{u_m}(\xi)
e^{2\pi i\xi
t}d\xi\right|^2\\&\leq\sum_m2^{j}\int_{2^{j-1}\leq|\xi|\leq2^j}|\widehat{u_m}(\xi)|^2
d\xi=2^j\left(\sum_m\|\Delta_iu_m\|_{L^2_t}^2\right).\end{split}\]

We distinguish three cases$$(\mathrm{i})\ \ |l|^{2N+1}\leq2^q\ \ \ \
(\mathrm{ii})\ \ |l|^{2N}\leq2^q<|l|^{2N+1}\ \ \ \ (\mathrm{iii})\ \
2^q<|l|^{2N}.$$

Contribution of $(\mathrm{i})$. By (\ref{8}), we have the estimate
\[\begin{split}\left\|\sum_{m}e^{2\pi iQ(m,l)t}
(\Delta_{p}u_m\overline{\Delta_{q}v_{m+l}})(t)\right\|_{L_t^2}&\leq\left\|(\sum_{m}|
\Delta_{p}u_m|^2)^{1/2}(\sum_{m}| \Delta_{p}v_{m+l}|^2)^{1/2}
\right\|_{L_t^2}\\&\leq2^{p/2}\left(\sum_m\|\Delta_pu_m\|_{L^2_t}^2\right)^{1/2}\left(\sum_m\|\Delta_qv_{m+l}
\|_{L^2_t}^2\right)^{1/2}.
\end{split}\]By Plancherel theorem we get
\begin{equation}\begin{split}\label{9}\left\|\sum_{m}e^{2\pi iQ(m,l)t}
(\Delta_{p}u_m\overline{\Delta_{q}v_{m+l}})(t)\right\|_{L_t^2}\lesssim&2^{p\frac{N}{4N+2}}2^{-q\frac{N+1}{4N+2}}
\left(\sum_m\int_{2^{p-1}\leq|\xi|\leq2^p}
(1+|\xi|)^{\frac{N+1}{2N+1}} |\widehat{u_m}(\xi)|^2d\xi
\right)^{1/2}\\&\times\left(\sum_m\int_{2^{q-1}\leq|\xi|\leq2^q}
(1+|\xi|)^{\frac{N+1}{2N+1}} |\widehat{v_m}(\xi)|^2d\xi
\right)^{1/2}.\end{split}\end{equation} Estimate (\ref{9}) implies
\[\begin{split}\sum_{q\geq p}\left\|\sum_{m}e^{2\pi iQ(m,l)t}
(\Delta_{p}u_m\overline{\Delta_{q}v_{m+l}})(t)\right\|_{L_t^2}\lesssim&\sum_{
2^p\geq|l|^{2N+1}}2^{-p\frac{1}{4N+2}}\left(\sum_m\int_{2^{p-1}\leq|\xi|\leq2^p}
(1+|\xi|)^{\frac{N+1}{2N+1}} |\widehat{u_m}(\xi)|^2d\xi
\right)^{1/2}\\&\times\left(\sum_mB_m^2\right)^{1/2}\\
\lesssim&|l|^{1/4}\left(\sum_{
2^p\geq|l|^{2N+1}}2^{-p\frac{1}{4N+2}}\sum_m\int_{2^{p-1}\leq|\xi|\leq2^p}
(1+|\xi|)^{\frac{N+1}{2N+1}} |\widehat{u_m}(\xi)|^2d\xi
\right)^{1/2}\\&\times\left(\sum_mB_m^2\right)^{1/2}.\end{split}\]
The last step is followed by H\"older inequality. For $j\geq0$, let
positive integer $|l|=2^j$ and $p^{\prime}=p-j(2N+1)$, we can obtain
the estimate
\[\begin{split}\sum_{l\in\mathbb{Z}}\left(\sum_{q\geq
p}\left\|\sum_{m\in\mathbb{Z}}e^{2\pi iQ(m,l)t}
(\Delta_{p}u_m\overline{\Delta_{q}v_{m+l}})(t)\right\|_{L_t^2}\right)^2&\lesssim\sum_{l\in\mathbb{Z}}\sum_{
2^p\geq|l|^{2N+1}}\sum_m\int_{2^{p-1}\leq|\xi|\leq2^p}
(1+|\xi|)^{\frac{N+1}{2N+1}} |\widehat{u_m}(\xi)|^2d\xi
\\&\times|l|^{1/2}2^{-p\frac{1}{4N+2}}\sum_mB_m^2\\
&\thickapprox\sum_{j\in\mathbb{Z}}\sum_{
p^{\prime}\geq0}\sum_m\int_{2^{p^{\prime}+j(2N+1)-1}\leq|\xi|\leq2^
{p^{\prime}+j(2N+1)}} (1+|\xi|)^{\frac{N+1}{2N+1}}
\\&\times|\widehat{u_m}(\xi)|^2d\xi2^{-p^{\prime}\frac{1}{4N+2}}\sum_mB_m^2\\&=\sum_{
p^{\prime}\geq0}2^{-p^{\prime}\frac{1}{4N+2}}\sum_mA_m^2\sum_mB_m^2\lesssim\sum_mA_m^2\sum_mB_m^2.\end{split}\]
Similarly, we have
$$\sum_{l\in\mathbb{Z}}\left(\sum_{q\geq
p}\left\|\sum_{m\in\mathbb{Z}}e^{2\pi iQ(m,l)t}
(\Delta_{p}v_m\overline{\Delta_{q}u_{m+l}})(t)\right\|_{L_t^2}\right)^2\lesssim\sum_mA_m^2\sum_mB_m^2.$$

Contribution of $(\mathrm{ii})$. Consider the quantity
$$\left\|\sum_{m\in\mathbb{Z}}e^{2\pi iQ(m,l)t}
(\Delta_{p}u_m\overline{\Delta_{q}v_{m+l}})(t)\right\|_{L_t^2}.$$We
first assume $l\geq0$. For a large positive $K>0$, by Plancherel
theorem, we write

\begin{equation}\label{10}\left\|\sum_{|m|\leq K}e^{2\pi iQ(m,l)t}
(\Delta_{p}u_m\overline{\Delta_{q}v_{m+l}})(t)\right\|_{L_t^2}=\left\|\sum_{|m|\leq
K}
\widehat{(\Delta_{p}u_m}*\widehat{\overline{\Delta_{q}v_{m+l}})}(\xi-Q(m,l))\right\|_{L_{\xi}^2}.\end{equation}
Since $\mathrm{supp}\
\widehat{(\Delta_{p}u_m}*\widehat{\overline{\Delta_{q}v_{m+l}})}(\xi)\subset[-2^{q+1},2^{q+1}]$,
splitting the summation into $62^q/|l|^{2N}$ summations over
arithmetic progressions of increment $62^q/|l|^{2N}$, say
$\mathcal{M}_s$ for $s=1,\cdots,62^q/|l|^{2N}$. if $m_1,\
m_2\in\mathcal{M}_s$ then
$$m_1=n_162^q/|l|^{2N}+d,\ \ \ \ m_2=n_262^q/|l|^{2N}+d.$$where $d$,
$n_1$ and $n_2$ denote be integers and $n_1>n_2$, $d<62^q/|l|^{2N}$.
We write
$$(m+l)^k-m^k=\sum_{\alpha+\beta=k,\
\beta\geq1}a_{\alpha,\beta}m^{\alpha}l^{\beta}$$ and obviously
$a_{\alpha,\beta}\geq0$ for all indices $\alpha$, $\beta$. In this
case we have
\begin{equation}\begin{split}\label{11}|Q(m_1,l)-Q(m_2,l)|&=\left|\sum_{k=1}^{2N+1}\sum_{\alpha+\beta=k,\
\beta\geq1}c_ka_{\alpha,\beta}(m_1^{\alpha}-m_2^{\alpha})l^{\beta}\right|\\&\geq(m_1-m_2)l^{2N}\geq\max\{62^{q},l^{2N}\}.
\end{split}\end{equation}By
a orthogonality consideration and (\ref{10}) we write
\begin{equation}\begin{split}\label{12}
\left\|\sum_{m\in\mathcal{M}_s}
\widehat{(\Delta_{p}u_m}*\widehat{\overline{\Delta_{q}v_{m+l}})}(\xi-Q(m,l))\right\|_{L_{\xi}^2}^2&=
\sum_{m\in\mathcal{M}_s}
\left\|\widehat{(\Delta_{p}u_m}*\widehat{\overline{\Delta_{q}v_{m+l}})}(\xi-Q(m,l))\right\|_{L_{\xi}^2}^2\\&=
\sum_{m\in\mathcal{M}_s}
\left\|\widehat{(\Delta_{p}u_m}*\widehat{\overline{\Delta_{q}v_{m+l}})}(\xi)\right\|_{L_{\xi}^2}^2\\&=
\sum_{m\in\mathcal{M}_s}
\left\|\Delta_{p}u_m\overline{\Delta_{q}v_{m+l}}\right\|_{L_t^2}^2.
\end{split}\end{equation}
From Young inequality, (\ref{12}) and (\ref{8}) we get
\[\begin{split}\left\|\sum_{|m|\leq K}e^{2\pi iQ(m,l)t}
(\Delta_{p}u_m\overline{\Delta_{q}v_{m+l}})(t)\right\|_{L_t^2}^2&\leq6\frac{2^q}{|l|^{2N}}
\sum_{s=1}^{62^q/|l|^{2N}}\left\|\sum_{m\in\mathcal{M}_s}e^{2\pi
iQ(m,l)t}
(\Delta_{p}u_m\overline{\Delta_{q}v_{m+l}})(t)\right\|_{L_t^2}^2\\&=6\frac{2^q}{|l|^{2N}}
\sum_{s=1}^{62^q/|l|^{2N}}\sum_{m\in\mathcal{M}_s}\left\|
(\Delta_{p}u_m\overline{\Delta_{q}v_{m+l}})(t)\right\|_{L_t^2}^2\\&\leq6\frac{2^q}{|l|^{2N}}2^p\sum_m
\|\Delta_{p}u_m\|_{L_t^2}^2\|\Delta_{p}v_m\|_{L_t^2}^2\\&\leq6|l|^{-2N}2^{2Np/(2N+1)}2^{2Nq/(2N+1)}\sum_m\int_{2^{p-1}
\leq|\xi|\leq2^p}(1+|\xi|)^{\frac{N+1}{2N+1}}\\&\ \ \ \ \ \
\times|\widehat{u_m}(\xi)|^2d\xi\int_{2^{q-1}
\leq|\xi|\leq2^q}(1+|\xi|)^{\frac{N+1}{2N+1}}|\widehat{v_{m+l}}(\xi)|^2d\xi.\end{split}\]
Let $K\rightarrow\infty$, the estimate holds for the whole integer
set. The same estimate can be obtained if we interchange the role of
$u_m$ and $v_m$. As for the case $l<0$, by (\ref{7}) we can reduce
this case to $l\geq0$.

Therefore, we have
\[\begin{split}\sum_{l\in\mathbb{Z}}\left(\sum_{q\geq p,\
2^q<|l|^{2N+1}}\left\|\sum_{m\in\mathbb{Z}}e^{2\pi iQ(m,l)t}
(\Delta_{p}u_m\overline{\Delta_{q}v_{m+l}})(t)\right\|_{L_t^2}\right)^2&\leq\sum_m\sum_lA_m^2B_{m+l}^2\\&
=\sum_mA_m^2\sum_mB_m^2.\end{split}\]

Contribution of $(\mathrm{iii})$. We have known from (\ref{11}) that
the orthogonality seems more natural in this case. The arguments are
similar as the case of $(\mathrm{ii})$ and even simpler.
\end{proof}
The result is sharp. For example if we take
$$u_K(x,t)=\sum_{|n|\leq K}\int_{|\lambda|\leq K^{2N+1}}e^{2\pi
i(nx+\lambda t)}d\lambda.$$Then $\|u_K\|_{L^2_xL^2_t}\thickapprox
K^{N+1}$ and$$\|u_K\|_{L^4_xL^4_t}\thickapprox
K^{3/4}(K^{2N+1})^{3/4}\thickapprox K^{3(N+1)/2}.$$On the other
hand,\[\begin{split}\|u_K\|_{X^{0,\frac{N+1}{4N+2}}}&=\left(\sum_{|n|\leq
K}\int_{|\lambda|\leq
K^{2N+1}}(1+|\lambda+P(n)|)^{\frac{N+1}{2N+1}}|
\widehat{u_K}(n,\lambda)|^2d\lambda\right)^{1/2}\\&\thickapprox
K^{3(N+1)/2}\thickapprox\|u_K\|_{L^4_xL^4_t}.\end{split}\]

\end{document}